\documentclass[12pt,twoside]{article}
\usepackage[hmargin=0.8in,vmargin=0.9in]{geometry}
\usepackage{fancyhdr}
\usepackage{graphicx}
\usepackage{float}
\usepackage{subfigure}
\usepackage{amssymb}
\usepackage{amsmath}
\usepackage{amsfonts}
\usepackage{theorem}
\usepackage{mathrsfs}
\usepackage{mathtools}
\usepackage{bm}
\usepackage{color}
\usepackage{setspace}
\usepackage{epstopdf}
\usepackage{cite}

\usepackage{booktabs} 
\usepackage{array} 
\usepackage{paralist} 
\usepackage{verbatim} 
\usepackage{subfigure} 

{\theorembodyfont{\rmfamily}\newtheorem{remark}{Remark}[section]}

\newenvironment{acknowledgment}{{\flushleft \bf Acknowledgment:}}{}

\allowdisplaybreaks[1]

\numberwithin{equation}{section}
\numberwithin{figure}{section}
\numberwithin{table}{section}

\newcommand\eref[1]{(\ref{#1})}

\def\bE{{\mathbb{E}}}
\def\eps{\varepsilon}

\def\sp{\sigma^+}
\def\sm{\sigma^-}
\def\spm{\sigma^\pm}

\def\rrp{\rho^+}
\def\rrm{\rho^-}
\def\rrpm{\rho^\pm}

\begin{document}

\author{Alina Chertock\thanks{Department of Mathematics, North Carolina State University, Raleigh, NC 27695, USA;
{\tt chertock@math.ncsu.edu}},~{}
Alexander Kurganov\thanks{Mathematics Department, Tulane University, New Orleans, LA 70118, USA; {\tt kurganov@math.tulane.edu}},~{}
Anthony Polizzi\thanks{Mathematics Department, Tulane University, New Orleans, LA 70118, USA; {\tt apolizz@tulane.edu}},
~{}and
Ilya Timofeyev\thanks{Department of Mathematics, University of Houston, Houston, TX 77204, USA; {\tt ilya@math.uh.edu}}}
\title{Pedestrian Flow Models with Slowdown Interactions}
\date{}

\maketitle
\begin{abstract}
In this paper, we introduce and study one-dimensional models for the behavior of pedestrians in a narrow street or corridor. We begin at
the microscopic level by formulating a stochastic cellular automata model with explicit rules for pedestrians moving in two opposite
directions. Coarse-grained mesoscopic and macroscopic analogs are derived leading to the coupled system of PDEs
for the density of the pedestrian traffic. The obtained PDE system is of a mixed hyperbolic-elliptic type and therefore, we rigorously
derive higher-order nonlinear diffusive corrections for the macroscopic PDE model. We perform numerical experiments, which compare and
contrast the behavior of the microscopic stochastic model and the resulting coarse-grained PDEs for various parameter settings and
initial conditions. We also demonstrate that the nonlinear diffusion is essential for reproducing the behavior of the stochastic system in
the nonhyperbolic regime.  
\end{abstract}

\section{Introduction}\label{sec:intro}
In contrast with the considerable effort devoted to the modeling of vehicular traffic, modeling of pedestrian traffic received little
attention until fairly recently. In the past two decades, a variety of pedestrian traffic
and crowd dynamics models have been proposed, ranging from agent-based microscopic to macroscopic PDE models describing various phenomena
such as crowd behavior under panic, pedestrian planning, structural design, etc., see \cite{helbing92,HM95,BD08} and the review papers
\cite{bedo11,helbing2001} for examples and references.

Cellular automata (CA) has a long history of applications in different areas of science and engineering. In particular, the CA models have
been applied to the vehicular traffic to derive a coarse-grained PDE description with the look-ahead dynamics for the density of the car
traffic \cite{soka06}. Multilane and multiclass were also considered as extensions of the original model \cite{also08,duso07,ckp10,WW}. CA
microscopic models with empirical rules have also been used to simulate pedestrian movement \cite{blad01,fuis99,bksz01,guhu08}.

The main emphasis of the present work is on the connection between the microscopic CA models for the bi-directional pedestrian traffic and
its coarse-grained PDE analogs. In particular, we utilize the CA approach to formulate a one-dimensional (1-D) microscopic model for the
pedestrian motion in a narrow street of corridor and derive a corresponding PDE description for the density of the pedestrian flow. 
The present work is motivated by the recent pedestrian experiments and modeling discussed in \cite{degond11,ardm11}, where the authors
consider the pedestrian motion in a circle and model the crowd movement using a system of 1-D PDEs.

The major advantage of the CA formalism is that it allows for a systematic derivation of the coarse-grained dynamics. The main assumptions
about the traffic flow (vehicular or pedestrian) are build into the microscopic model. The derivation of the coarse-grained description
typically requires some simplifying assumptions about the statistical behavior of the microscopic model which can be numerically verified.

The main conceptual difference between the CA models of the vehicular traffic and the CA model studied in this paper is that the pedestrian
motion is bi-directional which leads to completely new paradigms. In particular, the coarse-grained description of the pedestrian traffic
becomes a system of conservation laws of a mixed hyperbolic-elliptic type. We demonstrate that this system may exhibit an unrealistic
nonhyperbolic behavior, depending on the magnitude of the density of the pedestrian traffic. To overcome this difficulty, we systematically
derive nonlinear diffusive corrections to the PDE model from the original microscopic description.

The paper is organized as follows. In Section \ref{sec:micro}, we introduce the microscopic CA model. In Sections \ref{sec:meso} and
\ref{sec:macro}, we discuss the derivation of the mesoscopic and macroscopic analogs for the density of the pedestrian traffic, as well as
derive the next-order nonlinear diffusive corrections for the macroscopic PDE model. In Section \ref{sec:simul}, we describe several sets 
of simulations, in which the microscopic and macroscopic PDE models are compared illustrating, in particular, the effect of diffusion for
the initial conditions leading to the regime of nonhyperbolicity of the inviscid macroscopic PDE model. Finally, in Appendix A, we provide a
brief description of the numerical method used to solve the derived systems of PDEs.

\section{Microscopic Pedestrian Model}\label{sec:micro}
\begin{figure}[htbp]
\centerline{\scalebox{0.8}{\includegraphics{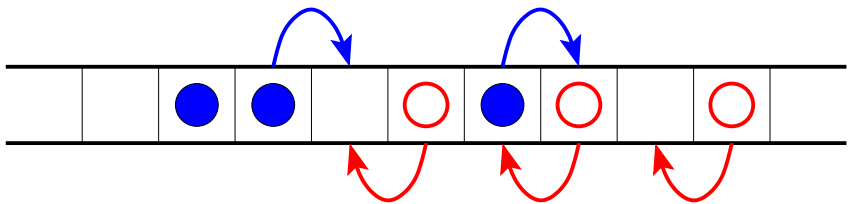}}}
\caption{\sf Schematic pedestrian configuration: Filled and empty circles represent pedestrians moving to the right and left, respectively.
Arrows represent possible pedestrian movements in this configuration.}
\label{fig1}
\end{figure}

To construct the microscopic model, we consider the time evolution of a 1-D lattice $\mathcal{L}$ with pedestrians moving into
two opposite directions. The process is schematically illustrated in Figure \ref{fig1}. 

The main difference with previous works on car traffic models (e.g., \cite{soka06}) is that the pedestrians can move into opposite
directions and we need to specify explicit rules when two such pedestrians meet. We represent pedestrians moving to the right and to the 
left by the variables $\sp_k(t)\equiv\sp(k,t)$ and $\sm_k(t)\equiv\sm(k,t)$, respectively, where $k$ is the position on the lattice and
$t$ denotes time: 
$$
\spm_k(t)=\left\{
\begin{aligned}
&1,&&\text{if at time $t$ cell $k$ is occupied by a pedestrian moving to the right (left)},\\
&0,&&\text{otherwise}. 
\end{aligned}\right.
$$
Moreover, we assume that pedestrians moving into two opposite directions can occupy the same cell (i.e., it is possible to have
$\sp_k(t)=\sm_k(t)=1$ for some $k$), but their velocities are reduced in that case. This mimics the slowdown due to side-stepping in more
realistic CA models \cite{blad01,fuis99} when pedestrians are avoiding the upcoming traffic. We also assume that pedestrians moving in the
same direction cannot occupy the same cell. For simplicity of the presentation, we omit the notation for the time-dependence of $\sp_k$ and
$\sm_k$ for the rest of the paper.

To complete the construction of the microscopic CA model, we consider explicit rules for the slowdown interaction which describe how the
velocity changes when pedestrians moving in opposite directions come in contact (i.e., occupy the same cell or two neighboring cells). In
particular, we explicitly prescribe velocities for four different pedestrian configurations in the cells neighboring to the right-moving
pedestrian with $\sp_k=1$ (assuming that $\sp_{k+1}=0$):
$$
\left\{\begin{aligned}
&c_0,&&\mbox{if $\sm_k=\sm_{k+1}=0$ (no left-moving pedestrians in cells $k$ or $k+1$)},\\
&c_1,&&\mbox{if $\sm_k=1, \sm_{k+1}=0$ (a left-moving pedestrian is in cell $k$)},\\
&c_2,&&\mbox{if $\sm_k=0, \sm_{k+1}=1$ (a left-moving pedestrian is in cell $k+1$)},\\
&c_3,&&\mbox{if $\sm_k=\sm_{k+1}=1$ (left-moving pedestrians in cells $k$ and $k+1$)}.
\end{aligned}\right.
$$
Notice that from the common sense considerations, the velocities should obey the following relationship: $c_3<c_2\approx c_1<c_0$.
The velocities of the left-moving pedestrian $\sm_k=1$ can be obtain in a similar manner.

The prescribed velocities are used to determine the probability of a pedestrian to move to the neighboring cell, that is, the probability
of transition $k\to k+1$ for $\sp_k$ during a small time interval $\Delta t$ is
\begin{equation}
\begin{aligned}
P^+_{k\to k+1}=\Delta t\Big[&c_0\sp_k(1-\sp_{k+1})(1-\sm_k)(1-\sm_{k+1})+c_1\sp_{k}(1-\sp_{k+1})\sm_k(1-\sm_{k+1})\\
+&c_2\sp_{k}(1-\sp_{k+1})(1-\sm_k)\sm_{k+1}+c_3\sp_{k}(1-\sp_{k+1})\sm_k\sm_{k+1}\Big],
\end{aligned}
\label{pplus}
\end{equation}
while the probability of transition $k\to k-1$ for $\sm_k$ is
\begin{equation}
\begin{aligned}
P^-_{k\to k-1}=\Delta t\Big[&c_0\sm_{k}(1-\sm_{k-1})(1-\sp_{k-1})(1-\sp_k)+c_1\sm_{k}(1-\sm_{k-1})(1-\sp_{k-1})\sp_k\\
+&c_2\sm_{k}(1-\sm_{k-1})\sp_{k-1}(1-\sp_k)+c_3\sm_{k}(1-\sm_{k-1})\sp_{k-1}\sp_k\Big].
\end{aligned}
\label{pminus}
\end{equation}

The lattice configurations $\bm{\sigma}^-_t:=\{\sm_k\}$ and $\bm{\sigma}^+_t:=\{\sp_k\}$ together constitute a continuous-time Markov chain
for $\bm{\sigma}_t:=\{\bm{\sigma}^-_t,\bm{\sigma}^+_t\}$. This model is easily simulated numerically using the Metropolis algorithm for
computing the expected values of $\sm$ and $\sp$. One can develop a Kinetic Monte-Carlo algorithm for these problems, but we found that the
Metropolis algorithm was quite efficient when the velocities $c_0, c_1, c_2, c_3$ are not very small.

Since $\bm{\sigma}_t$ is a continuous-time stochastic process, its generator is defined by
$$
A\psi=\lim_{\Delta t\to0}\frac{\bE\psi(\bm{\sigma}_{\Delta t})-\psi(\bm{\sigma}_0)}{\Delta t},
$$
where $\bm{\sigma}_0$ is the initial configuration, $\bm{\sigma}_{\Delta t}$ is the configuration at time $\Delta t$, $\psi$ is any test
function, and the expectation is taken over all possible transitions from $\bm{\sigma}_0$ to $\bm{\sigma}_{\Delta t}$. The generator of this
stochastic process can be computed as follows:
\begin{equation}
A\psi=\frac{1}{\Delta t}\sum_k\left\{P^+_{k\to k+1}\left[\psi(\bm{\sigma}^+_{k\leftrightarrow k+1},\bm{\sigma}^-_0)-\psi(\bm{\sigma}_0)
\right]+P^-_{k\to k-1}\left[\psi(\bm{\sigma}^+_0,\bm{\sigma}^-_{k\leftrightarrow k-1})-\psi(\bm{\sigma}_0)\right]\right\},
\label{gen}
\end{equation}
where $\bm{\sigma}^+_{k\leftrightarrow k+1}$ is the configuration obtained from $\bm{\sigma}^+_0$ by exchanging the values in cells $k$ and
$k+1$, and similarly, $\bm{\sigma}^-_{k\leftrightarrow k-1}$ is the configuration obtained from $\bm{\sigma}^-_0$ by exchanging the values
in cells $k$ and $k-1$. To derive the coarse-grain model, we need to compute the generator \eqref{gen} in two particular cases: For
$\psi(\bm{\sigma}^+,\bm{\sigma}^-)=\sigma^+_k$ and $\psi(\bm{\sigma}^+,\bm{\sigma}^-)=\sigma^-_k$, we obtain
\begin{equation}
A\sigma^+_k=\frac{P^+_{k-1\to k}-P^+_{k\to k+1}}{\Delta t}\quad\mbox{and}\quad A\sigma^-_k=\frac{P^-_{k+1\to k}-P^-_{k\to k-1}}{\Delta t},
\label{Asig}
\end{equation}
respectively.

\section{Mesoscopic Model}\label{sec:meso}
In this section, we use the microscopic CA model presented in Section \ref{sec:micro} to derive the mesoscopic model for the densities
$\rho_k^\pm:=\bE\spm_k$ (once again, to simplify notation, the time-dependence of $\rho_k^\pm$ is omitted throughout the paper).
To this end, we first recall that the generator $A$ satisfies the following property:
$$
\frac{d}{dt}\bE\psi=\bE A\psi,
$$
which can be applied to the test functions $\psi=\sp_k$ and $\psi=\sm_k$. This together with \eref{pplus}, \eref{pminus} and \eref{Asig}
results in the equations for the time-evolution of $\rho_k^\pm$:
\begin{equation}
\begin{aligned}
\frac{d\rho^+_k}{dt}=\bE\big[&c_0\sp_{k-1}(1-\sp_k)(1-\sm_{k-1})(1-\sm_k)-c_0\sp_k(1-\sp_{k+1})(1-\sm_k)(1-\sm_{k+1})\\
+&c_1\sp_{k-1}(1-\sp_k)\sm_{k-1}(1-\sm_k)-c_1\sp_{k}(1-\sp_{k+1})\sm_k (1-\sm_{k+1})\\
+&c_2\sp_{k-1}(1-\sp_k)(1-\sm_{k-1})\sm_k-c_2\sp_{k}(1-\sp_{k+1})(1-\sm_k)\sm_{k+1}\\
+&c_3\sp_{k-1}(1-\sp_k)\sm_{k-1}\sm_k-c_3\sp_{k}(1-\sp_{k+1})\sm_k\sm_{k+1}\big],\\
\frac{d\rho^-_k}{dt}=\bE\big[&c_0\sm_{k+1}(1-\sm_k)(1-\sp_k)(1-\sp_{k+1})-c_0\sm_{k}(1-\sm_{k-1})(1-\sp_{k-1})(1-\sp_k)\\
+&c_1\sm_{k+1}(1-\sm_k)(1-\sp_k)\sp_{k+1}-c_1\sm_{k}(1-\sm_{k-1})(1-\sp_{k-1})\sp_k\\
+&c_2\sm_{k+1}(1-\sm_k)\sp_k(1-\sp_{k+1})-c_2\sm_{k}(1-\sm_{k-1})\sp_{k-1}(1-\sp_k)\\
+&c_3\sm_{k+1}(1-\sm_k)\sp_k\sp_{k+1}-c_3\sm_{k}(1-\sm_{k-1})\sp_{k-1}\sp_k\big].
\end{aligned}
\label{micro1}
\end{equation}

The system \eqref{micro1} is exact, but not closed, since its right-hand side involves higher-order moments. The closure approximation
can be derived by assuming that the joint measure for $\sp$ and $\sm$ is approximately a product measure. This implies, for instance, that
adjacent cells are approximately independent and, in particular, $\bE[\spm_k\spm_{k+1}]\approx\bE\spm_k\bE\spm_{k+1}$. Thus, a closed system
of equations for $\rho_k^\pm$ can be obtained and the resulting mesoscopic model for the pedestrian density reads
\begin{equation}
\begin{aligned}
\frac{d\rrp_k}{dt}=&c_0\rrp_{k-1}(1-\rrp_k)(1-\rrm_{k-1})(1-\rrm_k)-c_0\rrp_{k}(1-\rrp_{k+1})(1-\rrm_k)(1-\rrm_{k+1})\\
+&c_1\rrp_{k-1}(1-\rrp_k)\rrm_{k-1}(1-\rrm_k)-c_1\rrp_{k}(1-\rrp_{k+1})\rrm_k (1-\rrm_{k+1})\\
+&c_2\rrp_{k-1}(1-\rrp_k)(1-\rrm_{k-1})\rrm_k-c_2\rrp_{k}(1-\rrp_{k+1})(1-\rrm_k) \rrm_{k+1}\\
+&c_3\rrp_{k-1}(1-\rrp_k)\rrm_{k-1}\rrm_k-c_3\rrp_{k}(1-\rrp_{k+1})\rrm_k\rrm_{k+1},\\
\frac{d\rrm_k}{dt}=&c_0\rrm_{k+1}(1-\rrm_k)(1-\rrp_k)(1-\rrp_{k+1})-c_0\rrm_{k}(1-\rrm_{k-1})(1-\rrp_{k-1})(1-\rrp_k)\\
+&c_1\rrm_{k+1}(1-\rrm_k)(1-\rrp_k)\rrp_{k+1}-c_1\rrm_{k}(1-\rrm_{k-1})(1-\rrp_{k-1})\rrp_k\\
+&c_2\rrm_{k+1}(1-\rrm_k)\rrp_k(1-\rrp_{k+1})-c_2\rrm_{k}(1-\rrm_{k-1})\rrp_{k-1}(1-\rrp_k)\\
+&c_3\rrm_{k+1}(1-\rrm_k)\rrp_k\rrp_{k+1}-c_3\rrm_{k}(1-\rrm_{k-1})\rrp_{k-1}\rrp_k.
\end{aligned}
\label{meso1}
\end{equation}
Note that the system \eqref{meso1} is defined on the same lattice $\mathcal{L}$ as the microscopic model.

\section{Macroscopic PDE Model}\label{sec:macro}
We now treat sites $k\in\mathcal{L}$ as cells with some fixed length $h>0$. Let $\Omega$ denote a subdomain of $\mathbb{R}$ corresponding
to the lattice $\mathcal{L}$, i.e., $\Omega=[0,L]$ (the number of cells, of course, depends on $h$). We consider a rescaling of time
$t\to ht$ and derive a coarse-grained PDE model in the limit as the cell size tends to zero and the number of cells tends to infinity. 

To this end, we rewrite the system \eqref{meso1} in the following flux form (taking the time rescaling into account):
\begin{equation}
\frac{d\rrp_k}{dt}=-\frac{F_{k,k+1}^+-F_{k-1,k}^+}{h},\qquad
\frac{d\rrm_k}{dt}=\frac{F_{k,k+1}^--F_{k-1,k}^-}{h},
\label{4.1}
\end{equation}
where
\begin{equation}
\begin{aligned}
&F_{k,k+1}^+=\rrp_k(1-\rrp_{k+1})\left[(1-\rrm_{k+1})\left(c_0(1-\rrm_k)+c_1\rrm_k\right)+
\rrm_{k+1}\left(c_2(1-\rrm_k)+c_3\rrm_k\right)\right],\\
&F_{k,k+1}^-=\rrm_{k+1}(1-\rrm_k)\left[(1-\rrp_k)\left(c_0(1-\rrp_{k+1})+c_1\rrp_{k+1}\right)+ 
\rrp_k\left(c_2(1-\rrp_{k+1})+c_3\rrp_{k+1}\right)\right].
\end{aligned}
\label{4.1a}
\end{equation}
Multiplying the above equations by $\varphi_k:=\varphi(kh)$, where $\varphi\in C^1_0(\overline{\Omega})$ is a test function, and using the
summation by parts property over $\Omega$, yields
\begin{equation}
\sum_k\varphi_k\frac{d\rrpm_k}{dt}=\pm\sum_kF_{k,k+1}^{\pm}\,\frac{\varphi_{k+1}-\varphi_k}{h}.
\label{4.2}
\end{equation}
Next, we multiply equation \eref{4.2} by $h$ and expand $\varphi_{k+1}$ in a Taylor series about $kh$ to obtain
\begin{equation}
\sum_k\varphi_k\frac{d\rrpm_k}{dt}h=\pm\sum_kF_{k,k+1}^{\pm}[\varphi_k'+\mathcal{O}(h)]h.
\label{4.3}
\end{equation}

We define pedestrian densities on $\Omega$ as follows. Again, using the notation $\rho^{\pm}$ (for convenience), define the function
$\rho^{\pm}(x,t)$ as a continuous piecewise linear interpolation (in the spatial variable) of $\rrpm_k(t)$ and take the limit as $h\to 0^+$.
Due to the boundeness of both $\rrpm$ and $\frac{d\rrpm_k}{dt}$ we obtain a weak formulation of the coarse-grained model:
\begin{equation}
\int\limits_{\Omega}\varphi(x)\frac{\partial}{\partial t}\rho^{\pm}(x,t)\,dx 
=\pm\int\limits_{\Omega}F^{\pm}(\rrp,\rrm)\varphi'(x)\,dx, 
\label{4.4}
\end{equation}
where $F^{\pm}(\rrp,\rrm)$ are defined as the corresponding limits of $F_{k,k+1}^{\pm}$, i.e.,
$$
F^+(\rrp,\rrm)=f(\rho^+)g(\rho^-),\quad F^-(\rrp,\rrm)=f(\rho^-)g(\rho^+),
$$ 
where
\begin{equation}
f(u)=u(1-u),\quad g(u)=(c_3-c_2-c_1+c_0)u^2+(c_2+c_1-2c_0)u+c_0.
\label{flux1}
\end{equation}
Since $\varphi$ is arbitrary, the integral equations \eref{4.4} can be written as the following system of PDEs:
\begin{equation}
\frac{\partial\rho^+}{\partial t}+\frac{\partial}{\partial x} \left[ f(\rho^+)g(\rho^-) \right]=0,\qquad
\frac{\partial\rho^-}{\partial t}-\frac{\partial}{\partial x} \left[ f(\rho^-)g(\rho^+) \right]=0.
\label{pde}
\end{equation}
\begin{remark}\label{rem41}
Note that the velocities $c_1$ and $c_2$ enter only as a sum into \eqref{flux1}. Therefore, it is not necessary to specify them separately.
\end{remark}
\begin{remark}
The coarse-grained system \eqref{pde} is only conditionally hyperbolic. Indeed, the Jacobian 
\begin{equation}
\left(
\begin{array}{rr}
f'(\rrp)g(\rrm)&f(\rrp)g'(\rrm)\\[0.8ex]-f(\rrm)g'(\rrp)&-f'(\rrm)g(\rrp)
\end{array}
\right)
\label{4.8}
\end{equation}
has real eigenvalues only if
\begin{equation}
\left[ f'(\rrm)g(\rrp) + f'(\rrp) g(\rrm) \right]^2 - 4 f(\rrm) f(\rrp) g'(\rrm) g'(\rrp) > 0.
\label{hypcondition}
\end{equation}
Therefore, for any particular choice of velocities $c_0, c_1, c_2$, and $c_3$ there is a region on nonhyperbolicity in the $(\rrm, \rrp)$
plane. From \eref{flux1}, one can see that the nonhyperbolicty can only manifest itself when pedestrians moving in two opposite directions
are both present in a particular location.
Two examples of nonhyperbolic regions are plotted in Figure \ref{fig:nonhyp}. The nonhyperbolic region described in \eqref{hypcondition} 
depends only on the ratio of velocities $c_1/c_0$, $c_2/c_0$, and $c_3/c_0$, but not on the particular value of $c_0$. The region of 
nonhyperbolicity becomes larger as the slowdown effect becomes more pronounces (i.e., as the ratios $c_1/c_0$, $c_2/c_0$, and $c_3/c_0$ 
become smaller). The loss of hyperbolicity may induce instabilities (as illustrated in Section \ref{sec:simul}), which are nonphysical 
and can be removed by introducing a nonlinear diffusive correction to the system (as we demonstrate in the next section). 
\begin{figure}[htbp]
\centerline{
\includegraphics[width=6.0cm,height=6.0cm]{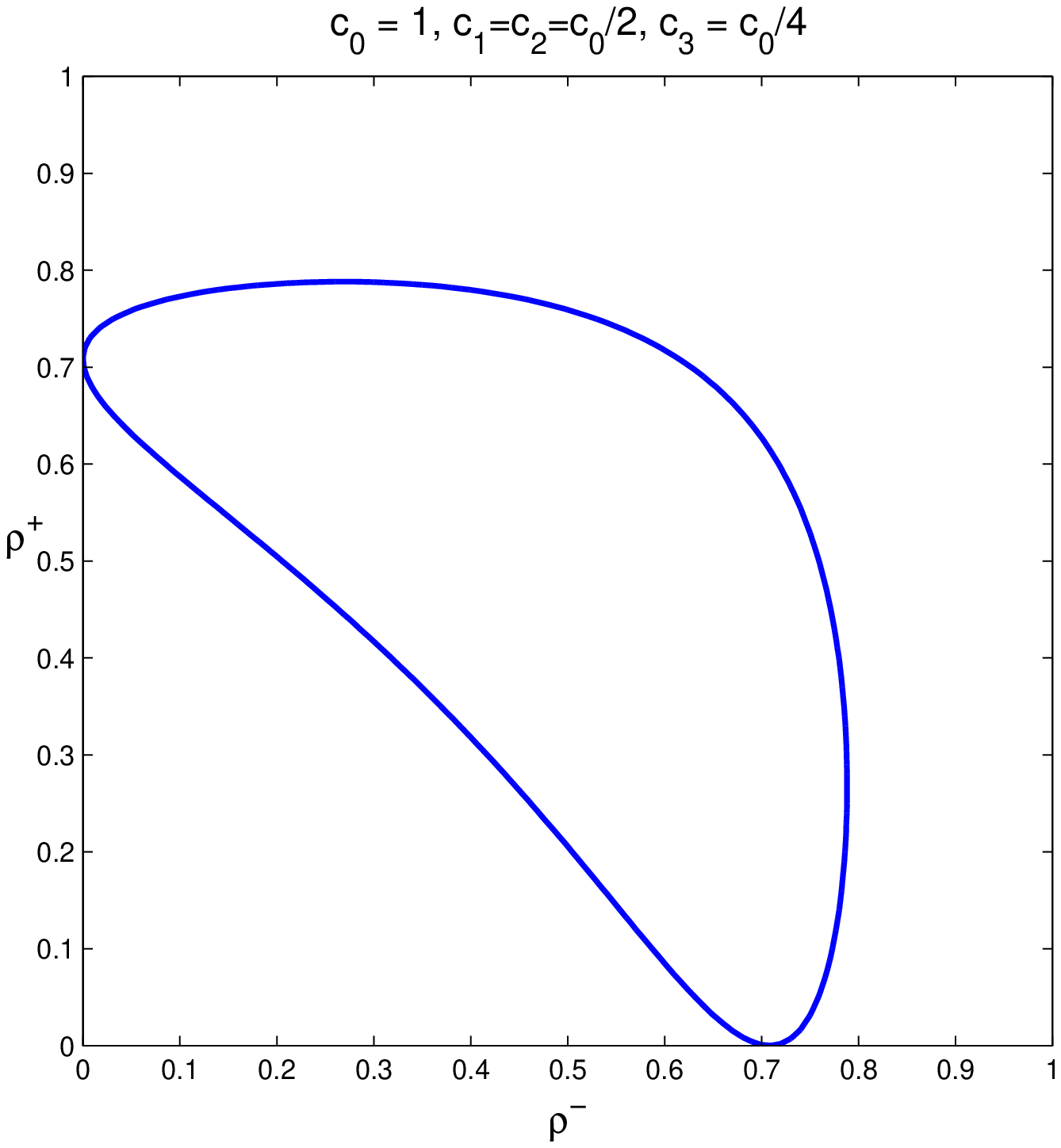}\hspace*{2cm}\includegraphics[width=6.0cm,height=6.0cm]{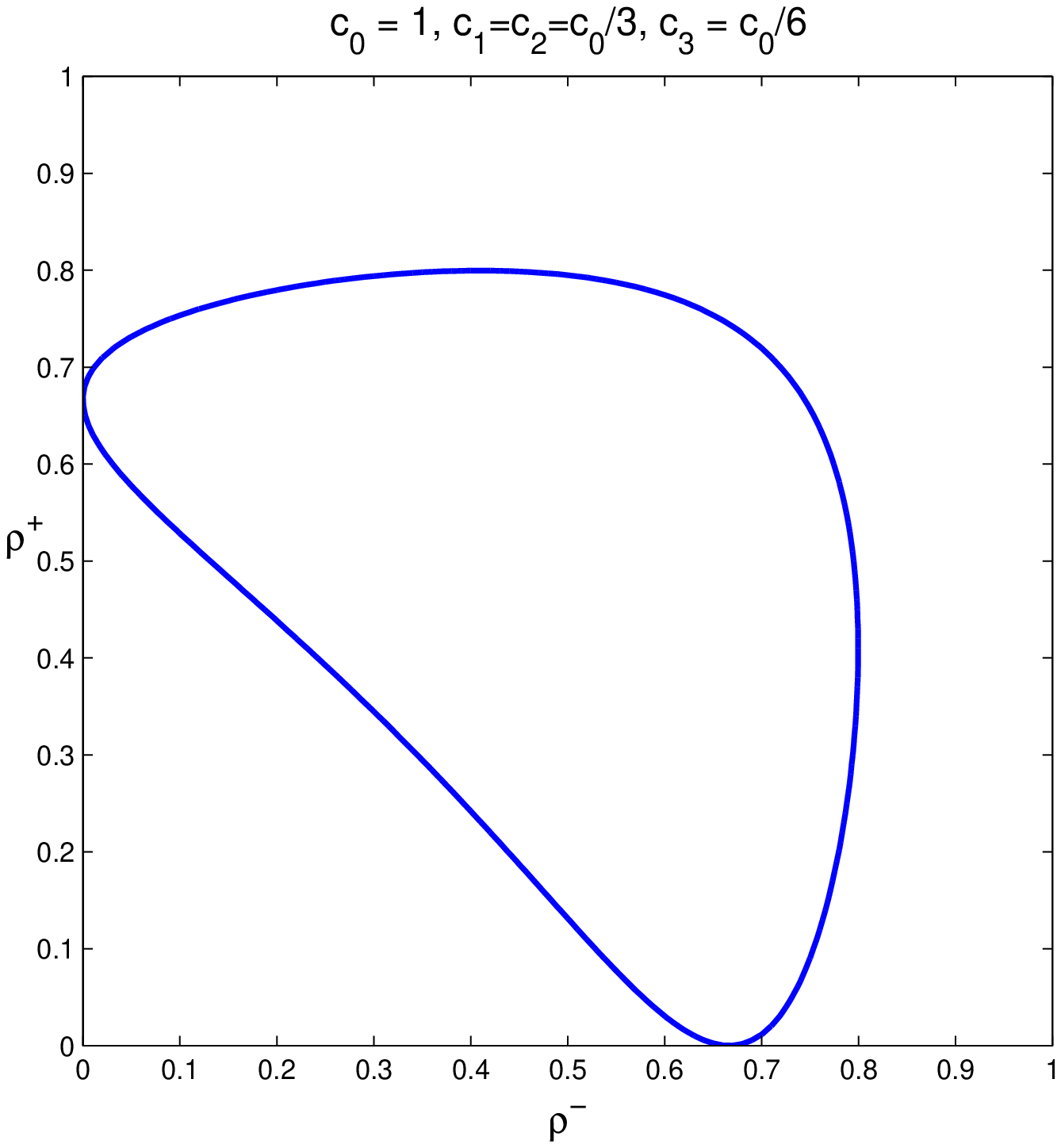}}
\caption{\sf Regions of nonhyperbolicity (inside) in the $\rrp$,$\rrm$ plane for the equation 
\eqref{pde} for two particular choices of velocities; left part $c_0=1$, $c_1=c_2=c_0/2$, $c_3=c_0/4$ and right part $c_0=1$, $c_1=c_2=c_0/3$, $c_3=c_0/6$.}
\label{fig:nonhyp}
\end{figure}
\end{remark}

\subsection{Diffusive Correction}\label{sec:diff}
The derivation of the coarse-grained PDE system \eref{pde}, \eref{flux1} can also be obtained by formally using the Taylor expansions
$$
\rrpm_{k\pm1}=\rrpm_k\pm h(\rrpm_k)'+\frac{h^2}{2}(\rrpm_k)''+\mathcal{O}(h^3)
$$ 
in the flux formulation \eref{4.1}, \eref{4.1a} followed by passing to the limit as $h\to 0^+$. 

Alternatively, keeping $h$ fixed and neglecting the $\mathcal{O}(h^3)$ terms, leads to the following second-order PDE system, which 
contains nonlinear diffusion terms: 
\begin{eqnarray*}
\frac{\partial\rho^+}{\partial t}+\frac{\partial}{\partial x}\left[f(\rho^+)g(\rho^-)\right]=h\Big[\frac{c_0}{2}\rrp_{xx}+ 
\left(c_1-c_0+(c_3-c_2-c_1+c_0)\rrm+(c_2-c_1)\rrp\right)\rrm_x\rrp_x\!\!\!&+&\\
\frac{1}{2}(c_1-c_2)\rrp (1-\rrp) \rrm_{xx}+
\frac{1}{2}\left((c_1 + c_2 - 2c_0)\rrm + (c_3-c_2-c_1+c_0)(\rrm)^2 \right) \rrp_{xx}\Big],\!\!\!&&\\
\frac{\partial\rho^-}{\partial t}-\frac{\partial}{\partial x}\left[f(\rho^-)g(\rho^+)\right] =
h\Big[\frac{c_0}{2} \rrm_{xx} + 
\left(c_1-c_0+(c_3-c_3-c_1+c_0)\rrp+(c_2-c_1)\rrm \right)\rrp_x \rrm_x\!\!\! &+& \\
\frac{1}{2}(c_1-c_2)\rrm (1-\rrm)\rrp_{xx} +
\frac{1}{2}\left((c_2+c_1-2c_0)\rrp+(c_3-c_2-c_1+c_0)(\rrp)^2\right)\rrm_{xx}\Big].\!\!\! &&
\end{eqnarray*}
The diffusion in the above system can be simplified considerably for certain conditions on the velocities. In particular, assuming that
\begin{equation}
c_1 = c_2,
\label{dc1}
\end{equation}
and replacing $h$ (a fixed cell size in the mesoscopic model) with a small parameter $\eps$, the diffusive correction to the original
coarse-grained PDE model \eref{pde} becomes
\begin{equation}
\begin{aligned}
\frac{\partial\rho^+}{\partial t}+\frac{\partial}{\partial x}\left[f(\rho^+)g(\rho^-)\right]&
=\frac{\eps c_0}{2}\left[\left((1-\rho^-)^2+2\alpha_1\rho^-(1-\rho^-)+\alpha_3(\rho^-)^2\right)\rrp_x\right]_x,\\[1.0ex]
\frac{\partial\rho^-}{\partial t}-\frac{\partial}{\partial x}\left[f(\rho^-)g(\rho^+)\right]&
=\frac{\eps c_0}{2}\left[\left((1-\rho^+)^2+2\alpha_1\rho^+(1-\rho^+)+\alpha_3(\rho^+)^2\right)\rrm_x\right]_x,
\end{aligned}
\label{pde2}
\end{equation}
where $\alpha_1=c_1/c_0$ and $\alpha_3=c_3/c_0$.
\begin{remark}
Notice that the assumption \eref{dc1} is rather mild since the velocities $c_1$ and $c_2$ only enter as a sum into the fluxes (see Remark
\ref{rem41}). We would also like to stress that the coefficients of the nonlinear diffusion are positive provided both $\rrp$ and $\rrm$ are
between 0 and 1.
\end{remark}
\begin{remark}
The nonlinear diffusion in \eqref{pde2} reflects the presence of pedestrians moving in the opposite direction. For instance, the diffusion
coefficient in the first, $\rrp$-equation depends on $\rrm$. It is instructive to consider the limiting cases: If $\rrm=0$ (i.e., no
pedestrians moving to the left are present), the diffusion reduces to the usual linear diffusion $0.5\eps c_0\rrp_{xx}$, while if $\rrm=1$, 
then the diffusion becomes $0.5\eps c_0\alpha_3\rrp_{xx}=0.5\eps c_3\rrp_{xx}$. The latter is, again, a linear diffusion, but with
a smaller coefficient (since $c_3<c_0$) reflecting a high density presence of the pedestrians moving in the opposite direction.
\end{remark}
\begin{remark}
It should be observed that the size of the viscosity coefficient $\eps$ needs to be established experimentally for each particular problem
at hand. In our numerical examples reported in Sections \ref{mixed} and \ref{sec53}, the value of $\eps$ is chosen empirically. 
\end{remark}

\section{Simulations}\label{sec:simul}
In this section, we present several sets of numerical experiments comparing and contrasting the propagation of pedestrian density in
ensemble simulations of the microscopic stochastic model and the corresponding numerical solutions of the macroscopic PDE models (with and
without diffusion). The mesoscopic model \eqref{meso1} can be considered as a conservative first-order finite-difference discretization of
the PDE models and the results of the mesoscopic simulations (not shown in the paper), performed on a sufficiently fine grid, typically
agree very well with the PDE results. 

In all of the PDE simulations below, we implement a semi-discrete second-order central-upwind scheme from \cite{KLin,KNP} with the minmod 
parameter $\theta=1$, CFL number equal to 0.5. The scheme is briefly described in Appendix \ref{app}. It should be observed that the PDE 
simulations are performed on a much coarser mesh than the corresponding microscopic ones. We use periodic boundary conditions in all of the
presented numerical examples. 

In the first set of numerical experiments (Section \ref{redlight}), the initial conditions resemble the ``red light'' situation when a
group of tightly packed pedestrians with density $1$ is released at time $t=0$. In the second set (Section \ref{mixed}), we mimic the
pedestrian movements starting with fully mixed initial conditions sampled from a particular piecewise constant density. Finally, in Section 
\ref{sec53}, we illustrate the stabilizing effect of the diffusive corrections in the nonhyperbolic regime.

\subsection{``Red Light'' Initial Conditions}\label{redlight}
We first consider the microscopic CA model with the following initial conditions:
\begin{equation}
\sp(k,0) = 
\begin{cases}
1,&n_1\le k\le n_2, \\
0,&\mbox{otherwise},
\end{cases}\qquad
\sm(k,0) = 
\begin{cases}
1,&N-n_2\le k\le N - n_1,\\
0,&\mbox{otherwise},
\end{cases}
\label{red:light:ic}
\end{equation}
with $n_2-n_1 \ll N$. These initial conditions correspond to two (relatively small) groups of pedestrians standing still and starting to
move toward each other at time $t=0$. The velocities are taken as
\begin{equation}
c_0=0.8m/s,~c_1=c_2 = c_0 /a,~c_3 = c_0 / (2a),
\label{vel1}
\end{equation}
where the parameter $a$ describes the strength of the slowdown interactions, for which we consider two regimes with either $a=2$
or $a=3$. The parameters in the CA model are chosen so that the two groups are away from the boundary and, thus, boundary conditions do not
affect the interaction:
\begin{equation}
N=1400,~n_1=301,~n_2=340,~h=0.2m,~\Delta t=0.01s,~MC=5000,
\label{red:light:params}
\end{equation}
%
%
where $MC$ is the number of Monte-Carlo simulations.

The macroscopic simulations of the purely convective PDE system \eref{pde} are performed on the computational domain $[0,280]$ with the
mesh size $\Delta x=0.8$ and initial data corresponding to \eref{red:light:ic}:
\begin{equation}
\rrp(x,0) = 
\begin{cases}
1,&60<x<68, \\
0,&\mbox{otherwise},
\end{cases}\qquad
\rrm(x,0) = 
\begin{cases}
1,&212<x<220,\\
0,&\mbox{otherwise},
\end{cases}
\label{ic1}
\end{equation} 
with the velocities $c_0, c_1, c_2$ and $c_3$ same as in \eref{vel1}. 

\begin{figure}[htbp]
\centerline{\includegraphics[width=12cm]{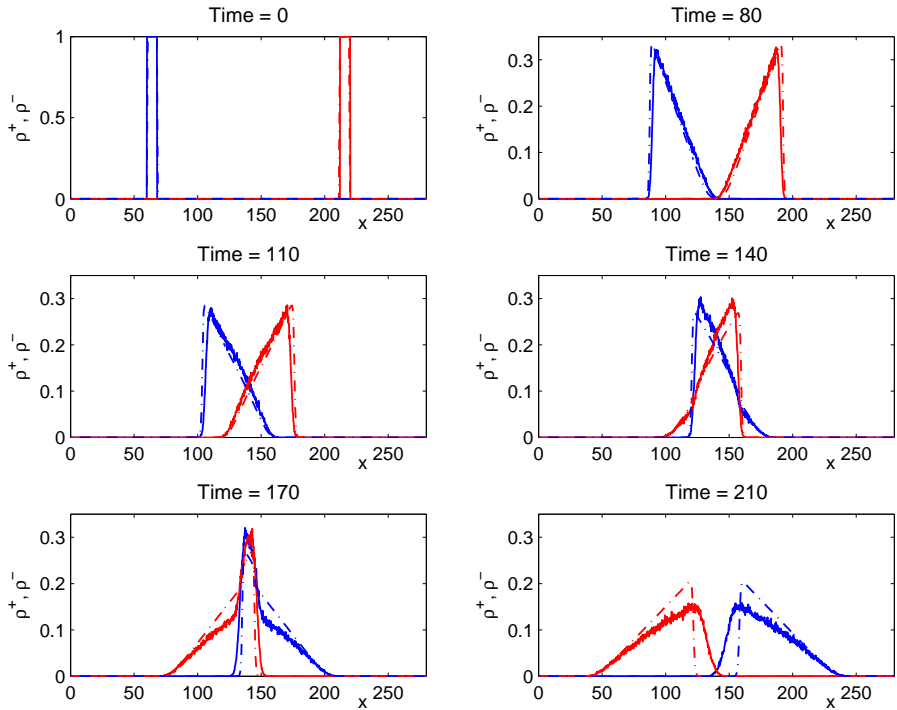}}
\caption{\sf Comparison of the density profiles computed by the CA (solid line) and PDE (dash-dotted line) models with $a=2$ and the 
``red light'' initial conditions in \eref{red:light:ic} and \eref{ic1}. The right- ($\rho^+$) and left-moving ($\rho^-$) pedestrians are
represented by the blue and red colors, respectively.}
\label{fig:alpha2}
\end{figure}

\begin{figure}[htbp]
\centerline{\includegraphics[width=12cm]{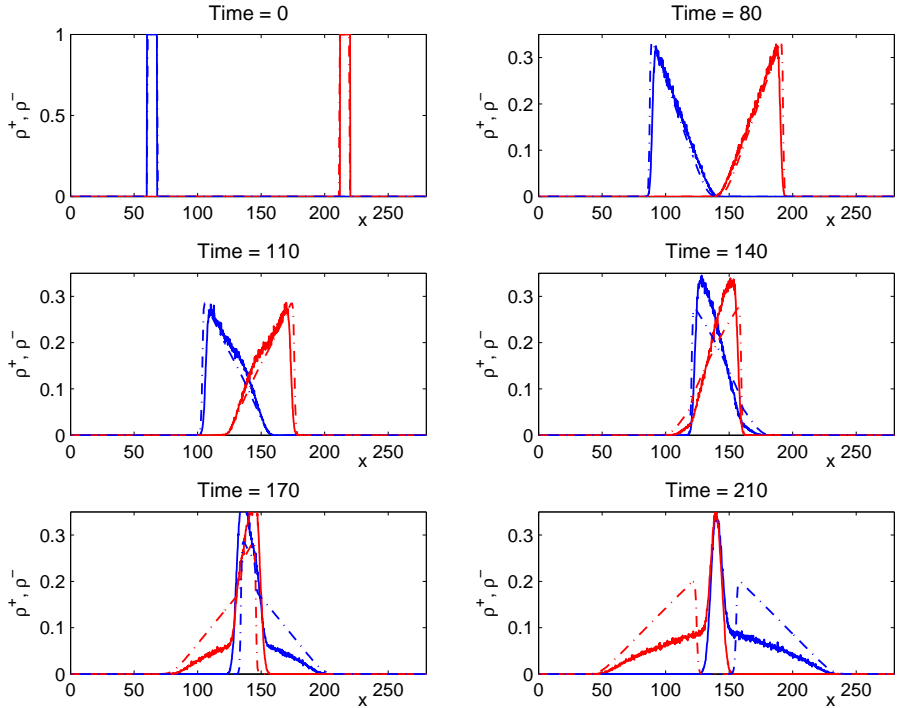}}
\caption{\sf The same as in Figure \ref{fig:alpha2} but with $a=3$.} 
\label{fig:alpha3}
\end{figure}

The initial settings \eref{red:light:ic}--\eref{red:light:params} and \eref{ic1} correspond to the right- and left-moving
groups, which initially do not overlap and their dynamics is equivalent to moving cars. After some time, the two groups begin to interact
and we study how well the macroscopic PDE model reproduces these interactions. Figures \ref{fig:alpha2} and \ref{fig:alpha3} show
comparison of the density profiles in CA and PDE simulations. 

The macroscopic PDE model \eqref{pde} produces very good results in the $a=2$ regime. In particular, it reproduces the forward front at
times $t=80, 110$ and $140$ quite well. The macroscopic results are slightly different from the stochastic ones after the two groups
pass through each other (see times $t=170$ and $210$): The back front of moving pedestrians in simulations of the macroscopic model is
sharper compared with the simulations of the stochastic model at these times. This can be explained by the purely hyperbolic nature of the
macroscopic PDE for this range of the density values.

The $a=3$ regime corresponds to stronger slowdown interactions between the two groups of pedestrians. We observe a considerable discrepancy
between the CA and PDE models at times $t=140, 170$ and $210$. In particular, an effect similar to blocking occurs in the 
Monte-Carlo simulations of the stochastic model; this is especially evident at time $t=170$. The propagation of two blocked groups
resembles a slow diffusion over a barrier, whereas the simulations of the PDE model exhibit a faster propagating wave. Therefore, the PDE
model with $a=3$ seems to be outside of the range of validity of the closure approximations used in its derivation.

We note that the results obtained using the viscous PDE model \eref{pde2} are quite similar to the inviscid ones: The diffusion simply
smears the fronts but no substantially different phenomena have been observed.

\subsection{Fully Mixed Initial Conditions}\label{mixed}
In the second set of numerical experiments, we simulate pedestrian movement in a periodic domain, which is divided into 30 sectors with 15
cells in each sector (totally $N=450$ cells). We initialize the CA simulations with the total number of 70 pedestrians (with 35 pedestrians
moving in each direction). The initial numbers of right- and left-moving pedestrians in each sector are specified, while their distribution
inside each sector is purely random (uniform). In dimensional units, each cell is $0.466m$ long, each sector is $7m$ long, and the domain is
$210m$ long. Other parameters in the CA simulations are chosen to be $\Delta t=0.005$ and $MC=3000$.

%
%
%
%

The corresponding initial densities for the PDE models \eref{pde} and \eref{pde2} are
\begin{equation}
\rrpm(x,0)=\frac{n^{\pm}_i}{15}~~\mbox{for}~~\frac{i-1}{30}L<x<\frac{i}{30}L,\quad i=1,\ldots,30,
\label{ic2}
\end{equation}
where $n^+_i$ and $n^-_i$ are numbers of right- and left-moving pedestrians in the $i$th sector, respectively. The size of the computational 
domain is $L=210$ and $\Delta x=1$. In the simulations of the viscous PDE \eref{pde2}, we take the diffusion coefficient $\eps=0.5$.

As in Section \ref{redlight}, we take the velocities 
\begin{equation}
c_0=1{m/s},\quad c_1=c_2=c_0/a,\quad c_3=c_0/(2a),
\label{weak} 
\end{equation}
and perform two sets of numerical simulations with velocities with $a=2$ and $a=3$.

\begin{figure}[htpb]
\centerline{\includegraphics[width=16cm,height=7.5cm]{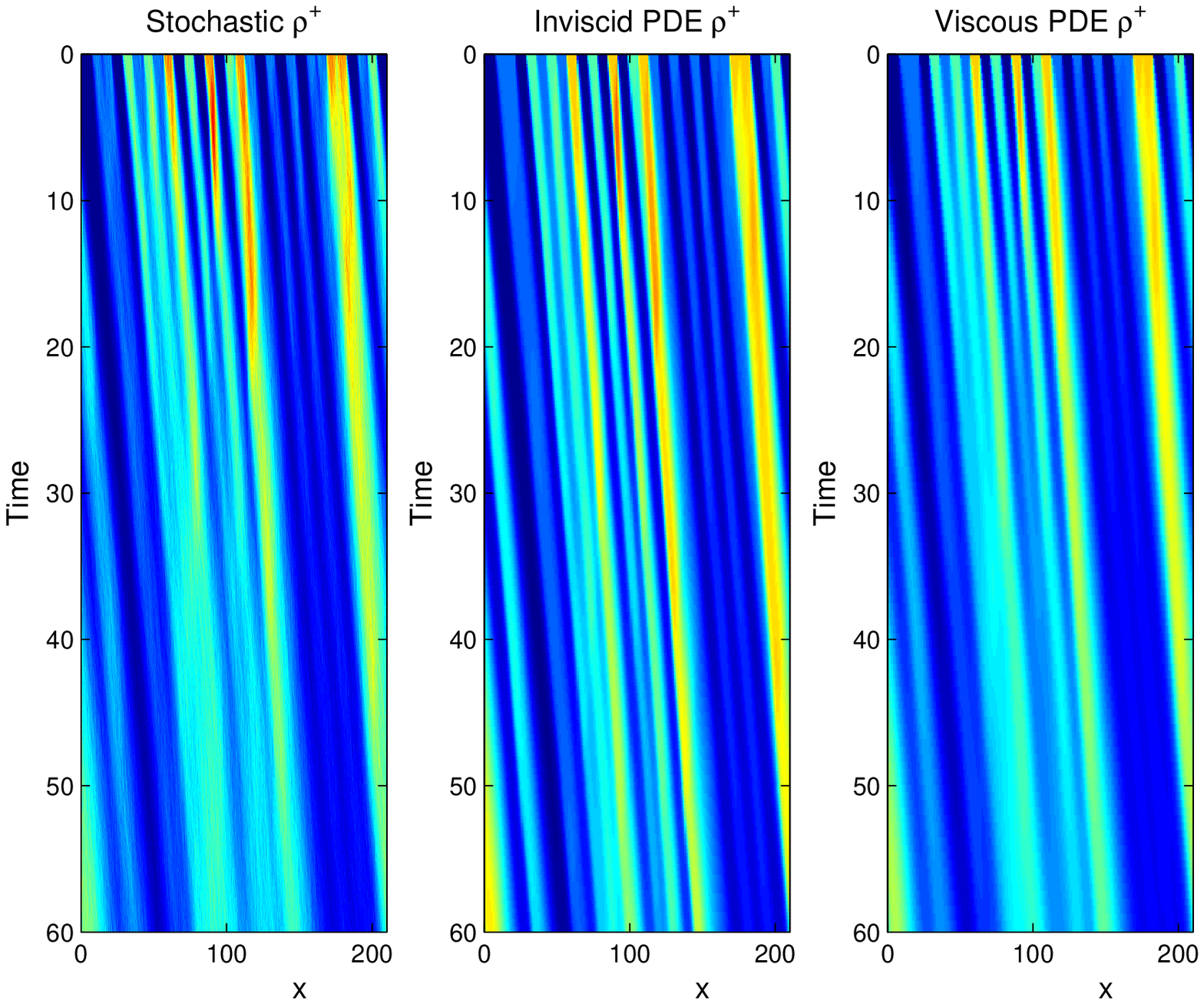}}
\caption{\sf Dynamics of right-moving pedestrians in both the CA and PDE models with the fully mixed initial conditions \eref{ic2} and
velocities \eref{weak} with $a=2$.}
\label{fm1}
\end{figure}

\begin{figure}[ht!]
\centerline{\includegraphics[width=16cm,height=11.0cm]{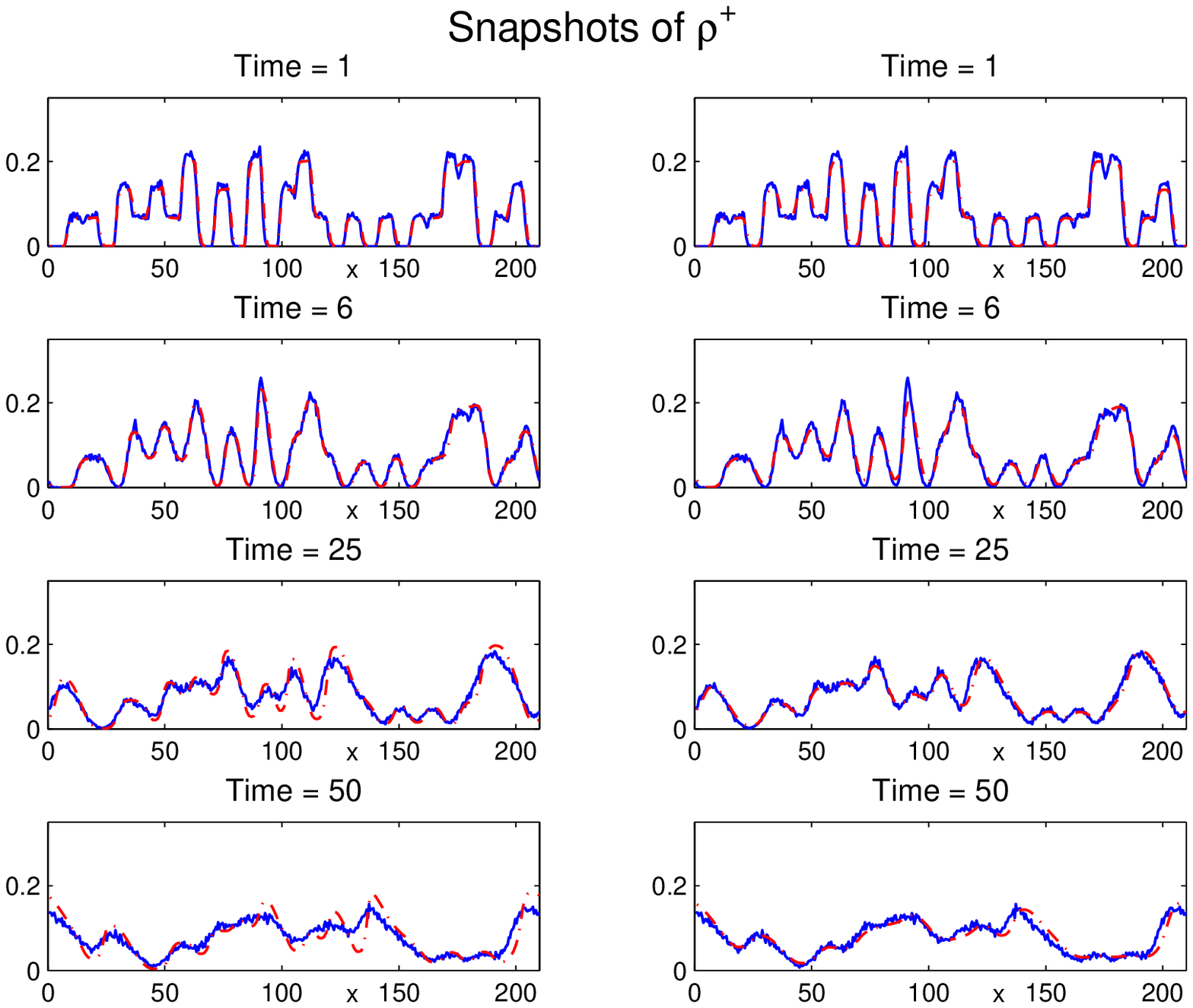}}
\caption{\sf Time snapshots of the solutions from Figure \ref{fm1}. The solid blue line represents the CA solution, while the dash-dotted
red line represents the inviscid (left) and viscous (right) PDE solutions.}
\label{fm2}
\end{figure}

\begin{figure}[htpb]
\centerline{\includegraphics[width=16cm,height=7.5cm]{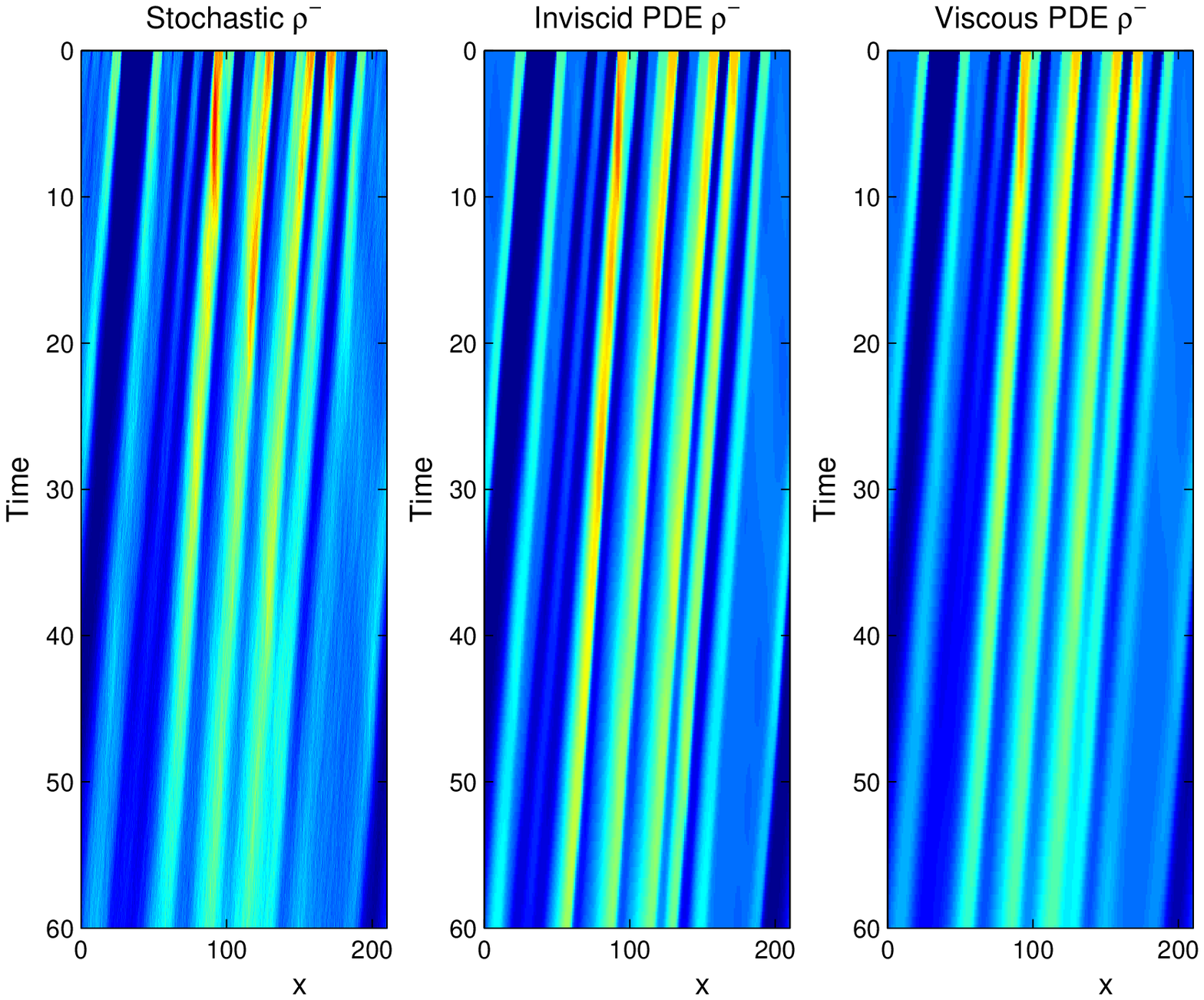}}
\caption{\sf The same as in Figure \ref{fm1} but for left-moving pedestrians.}
\label{fm3}
\end{figure}

\begin{figure}[htbp]
\centerline{\includegraphics[width=16cm,height=11.0cm]{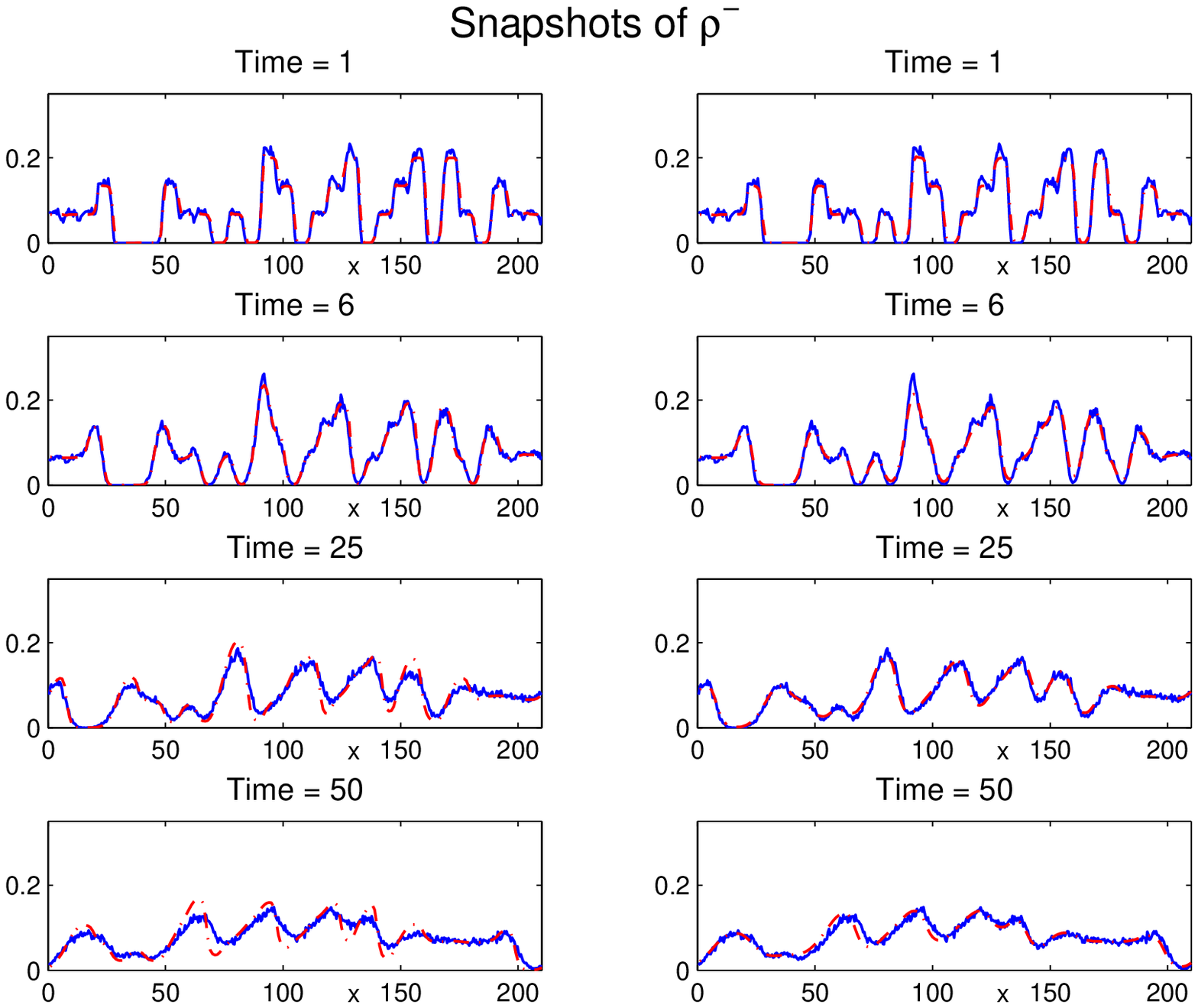}}
\caption{\sf Time snapshots of the solution from Figure \ref{fm3}. The solid blue line represents the CA solution, while the dash-dotted 
red line represents the inviscid (left) and viscous (right) PDE solutions.}
\label{fm4}
\end{figure}

\begin{figure}[htbp]
\centerline{\includegraphics[width=16cm,height=7.5cm]{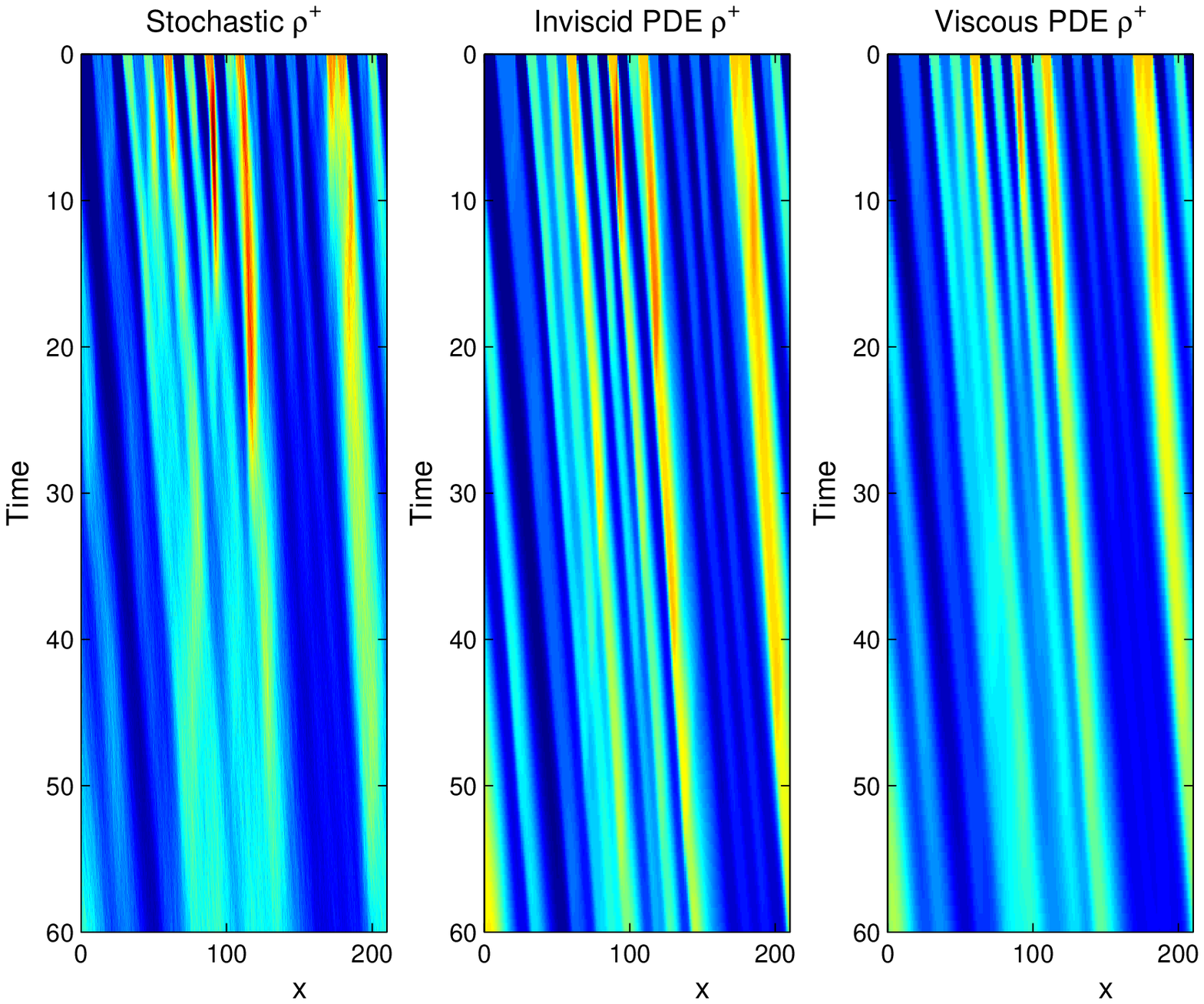}}
\caption{\sf Dynamics of right-moving pedestrians in both the CA and PDE models with the fully mixed initial conditions \eref{ic2} and
velocities \eref{weak} with $a=3$.}
\label{fm5}
\end{figure}

\begin{figure}[htbp]
\centerline{\includegraphics[width=16cm,height=11.0cm]{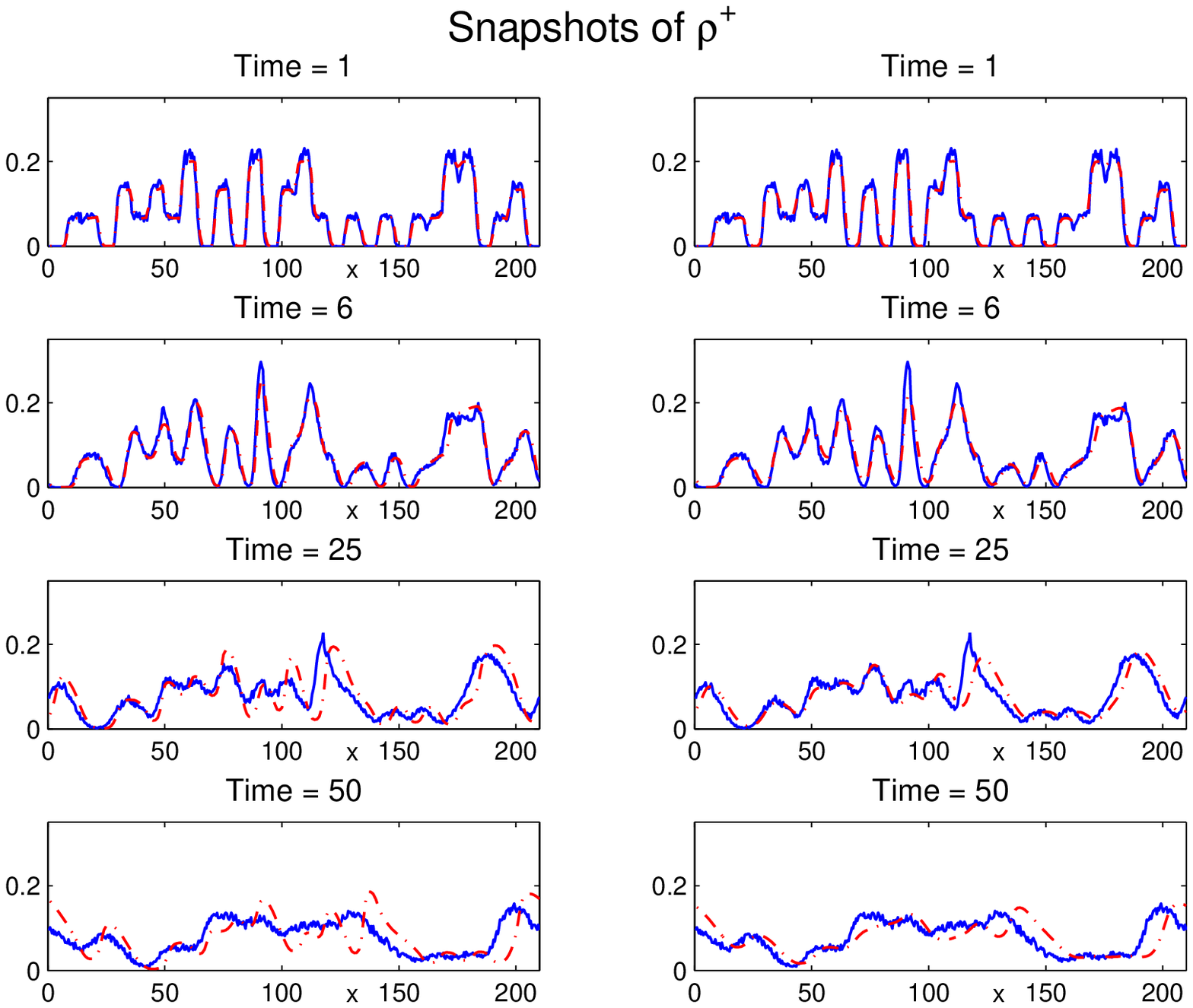}}
\caption{\sf Time snapshots of the solutions from Figure \ref{fm5}. The solid blue line represents the CA solution, while the dash-dotted 
red line represents the inviscid (left) and viscous (right) PDE solutions.}
\label{fm6}
\end{figure}

\begin{figure}[htbp]
\centerline{\includegraphics[width=16cm,height=7.5cm]{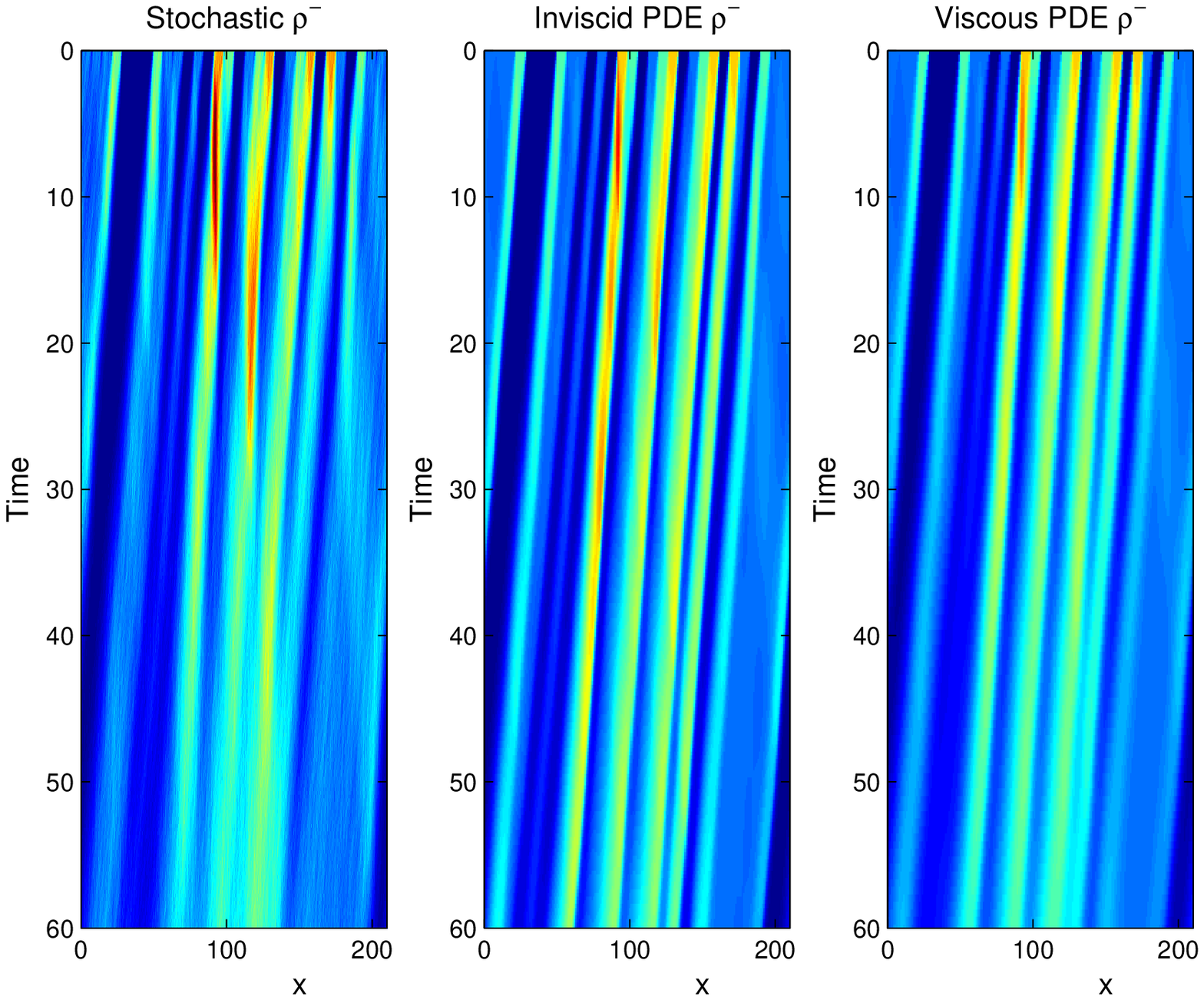}}
\caption{\sf The same as in Figure \ref{fm5} but for left-moving pedestrians.}
\label{fm7}
\end{figure}

\begin{figure}[ht!]
\centerline{\includegraphics[width=16cm,height=11.0cm]{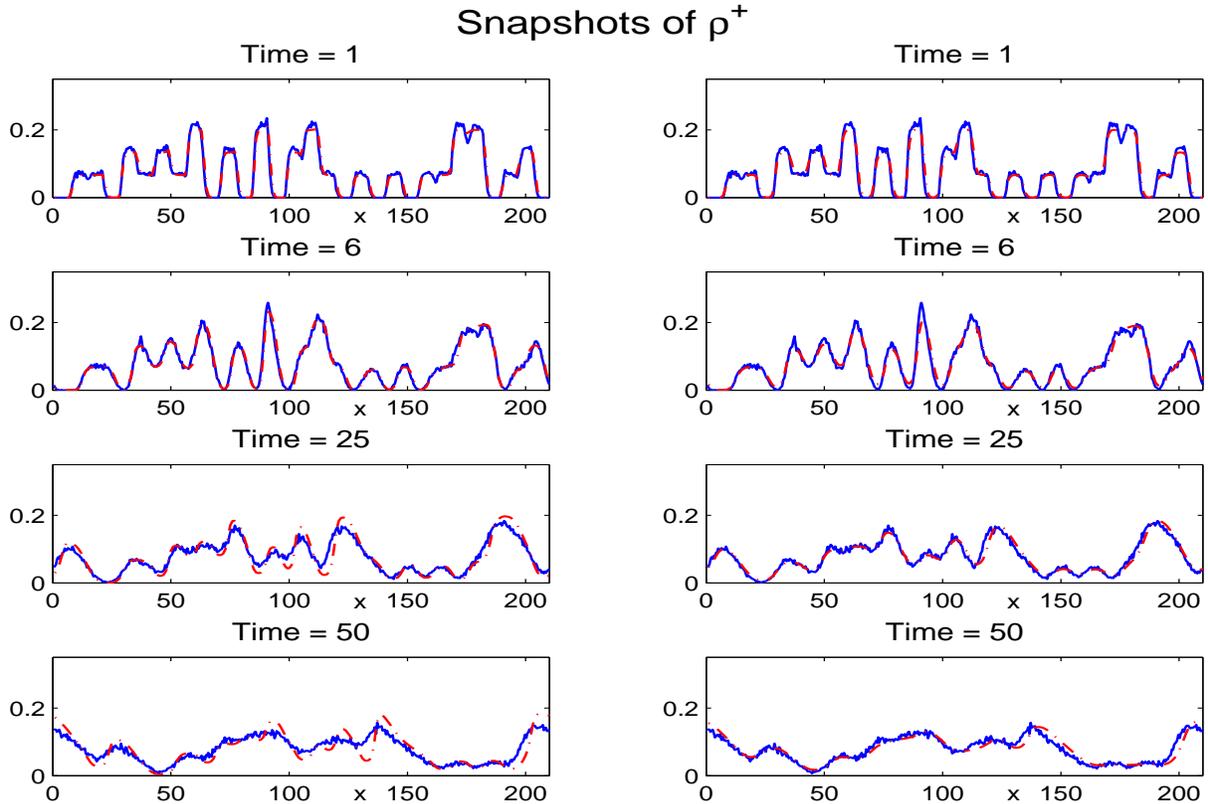}}
\caption{\sf Time snapshots of the solution from Figure \ref{fm6}. The solid blue line represents the CA solution, while the dash-dotted red
line represents the inviscid (left) and viscous (right) PDE solutions.}
\label{fm8}
\end{figure}

The obtained results are plotted in Figures \ref{fm1}--\ref{fm8}. In the $a=2$ case, the CA and PDE results are in a very good agreement, 
especially when the diffusion terms are included into the PDE models (Figures \ref{fm1}--\ref{fm4}). By comparing Figures \ref{fm1} and
\ref{fm3}, one can see that the invscid PDE model preserves initial pedestrian clusters longer than the CA or viscous PDE ones. The major
qualitative difference between the CA and viscous PDE models is a ``blocking'' phenomenon observed in the stochastic simulations. When the
slowdown interaction is stronger ($a=3$), the difference between the CA and PDE simulations are more pronounced: The stochastic ``blocking''
is more severe, while the PDE models develop both ``blocking'' and ``stop-and-go'' waves (see Figures \ref{fm5}--\ref{fm8}).

\subsection{Nonhyperbolic Regime}\label{sec53}
In both simulations described in Sections \ref{redlight} and \ref{mixed}, the magnitude of pedestrian densities $\rrp$ and $\rrm$ remain
smaller than the critical values for the nonhyperbolic behavior (recall that the hyperbolicity criterion is given by \eref{hypcondition}. 

To illustrate the nonhyperbolic behavior we choose the velocities to be \eref{weak} with $a=2$ and the initial density in the simulations 
of the coarse-grained PDE models to be
\begin{equation}
\rrp(x,0)=\begin{cases}
0.6,&140\le k\le 210,\\
0,&\mbox{otherwise},
\end{cases} \quad
\rrm(x,0)=\begin{cases}
0.6,&186.6\le x\le 233.3, \\
0, &\mbox{otherwise},
\end{cases}
\label{nonhic}
\end{equation}
that is, the initial data are in the nonhyperbolic regime. The computational domain is $[0,420]$ and the simulations are conducted on a
$\Delta x=420/1280$ grid.

The initial pedestrian distribution for the CA computations is sampled from the initial density \eref{nonhic}. The number of cells in 
the stochastic simulations is $N=900$, the cell size is $420/900\approx 0.4667m$, the time step is $\Delta t=0.005$ and $MC=3000$. 

\begin{figure}[ht!]
\centerline{\includegraphics[width=16cm,height=11.0cm]{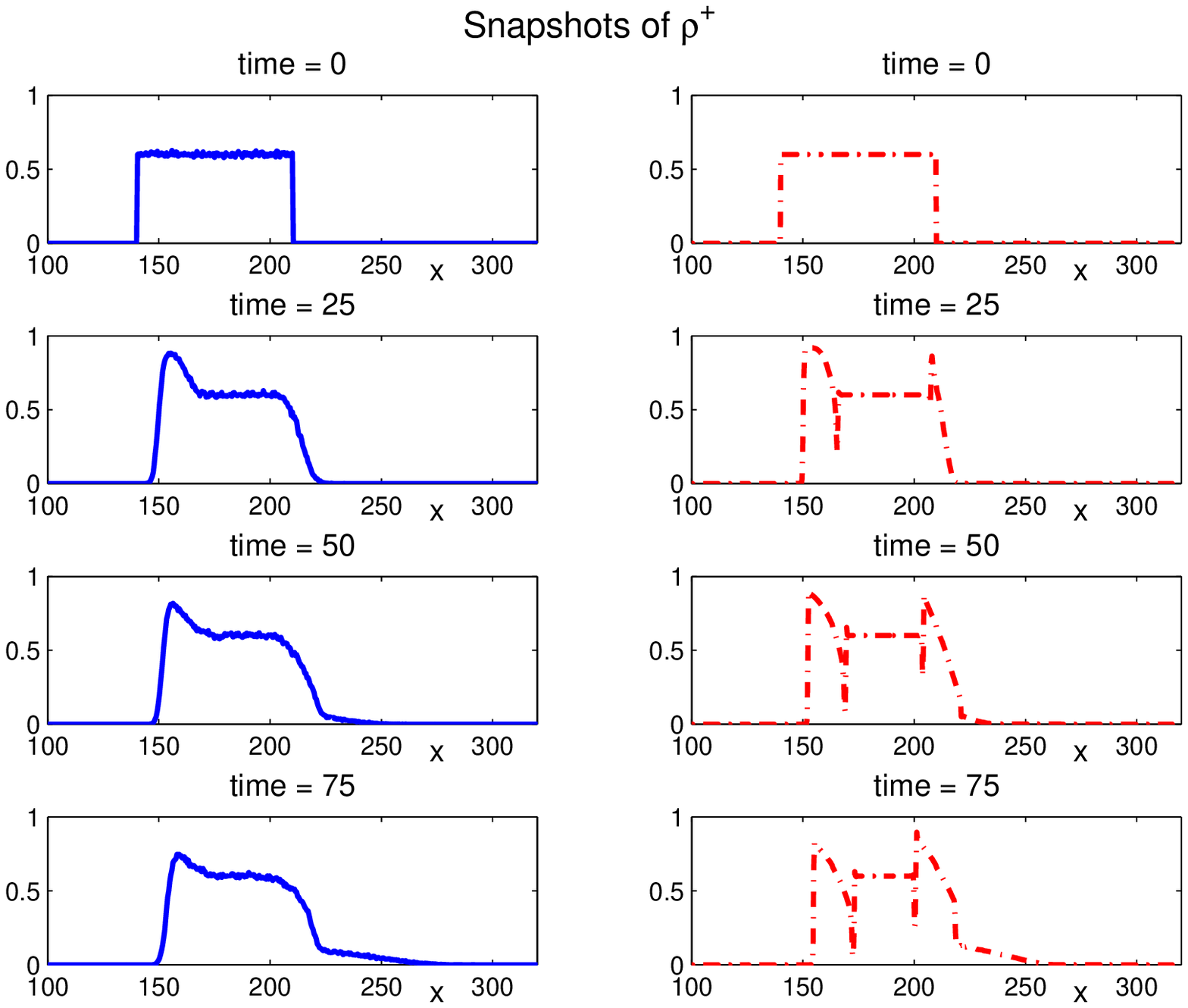}}
\caption{\sf Densities of the right-moving pedestrians computed by the CA (left) and inviscid PDE (right) models with the initial conditions
\eref{nonhic} and velocities \eref{weak} with $a=2$.}
\label{fig:nonh1}
\end{figure}

\begin{figure}[ht!]
\centerline{\includegraphics[width=16cm,height=11.0cm]{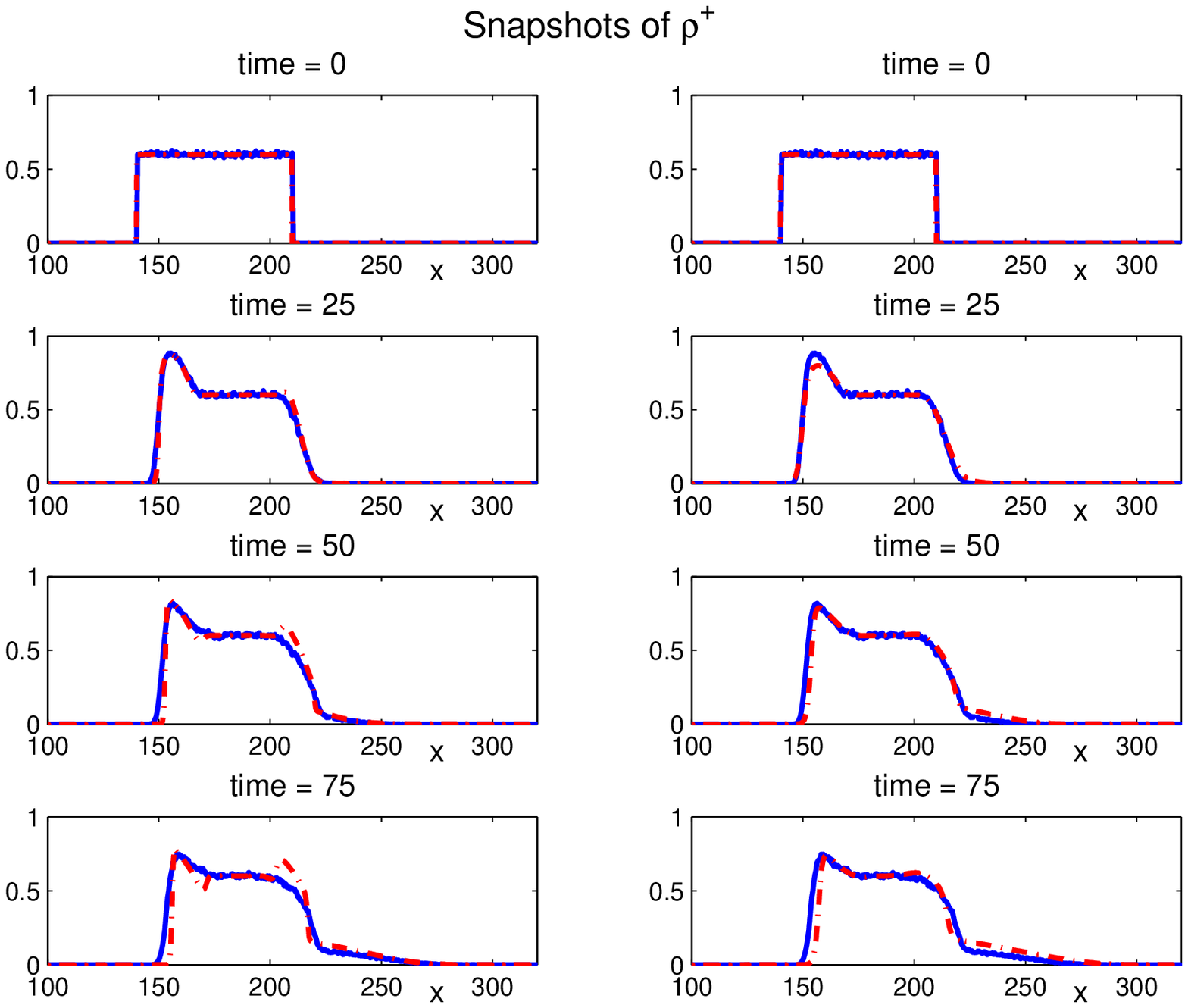}}
\caption{\sf Densities of the right-moving pedestrians computed by the CA and viscous PDE models with the initial conditions \eref{nonhic}
and velocities \eref{weak} with $a=2$. The solid blue line represents the CA solution, while the dash-dotted red line represents viscous 
PDE solutions with $\eps=0.5$ (left) and $\eps=1.5$ (right).}
\label{fig:nonh2}
\end{figure}

Here, we depict only $\rrp$ since $\rrm$ exhibits a similar behavior. Comparison between the stochastic and inviscid PDE simulations is 
shown in Figure \ref{fig:nonh1}. As one can see, the PDE solution develops spurious large magnitude oscillations, which demonstrates that
the inviscid system \eref{pde} is ill-posed. Figure \ref{fig:nonh2} shows $\rrp$ computed using the viscous PDE model \eref{pde2}. The
results suggest that the nonlinear diffusion present in \eref{pde2} stabilizes the PDE solution: The bigger value of $\eps=1.5$ fully
supresses spurious oscillations and leads to a very good agreement between the CA and PDE simulations.

\section{Conclusion}\label{conc}
In this paper, we have developed a PDE formulation for the pedestrian dynamics in a narrow street or corridor. The resulting 1-D system
of PDEs has been systematically derived from a microscopic stochastic CA model. The slowdown rules in the CA model need to be specified
explicitly by prescribing different velocities when pedestrians moving in the opposite directions are present. The condition of approximate
independence of neighboring cells is essential for the derivation of the coarse-grained PDE model and ``propagate'' into the macroscopic
PDE formulation resulting in a particular form of the flux. In addition, we have also demonstrated that the resulting PDE is only
conditionally hyperbolic. To stabilize the resulting PDE system we have rigorously derived the nonlinear diffusion correction using the
intermediate mesoscopic model. 

We have performed a number of numerical experiments comparing and contrasting the statistical behavior of the stochastic solution and its
deterministic PDE counterparts. Our numerical experiments indicate that in the case of weak slowdown interactions, the average behavior of
the stochastic CA model is reproduced quite well by the coarse-grained PDE model. On the other hand, the approximate independence assumption
does not seem to hold in the case of the stronger slowdown regime, which results in a considerable discrepancy between the stochastic and 
deterministic results. In addition, our numerical experiments also indicate that the nonlinear diffusive corrections may lead to a
significant improvement in reproducing the CA results using the PDE model. We would also like to point out that the PDE simulations are an
order of magnitude faster than the corresponding CA ones. Therefore, the coarse-grained PDE systems have considerable numerical advantages
in practical applications such as real-time predictions of the pedestrian traffic via filtering, investigation of pedestrian traffic flow
on networks, etc. These and other issues will be investigated in subsequent papers.

\begin{acknowledgment}
The first ideas for this work were discussed by the authors during their participation in the Program on Complex Networks organized by the
NSF funded Statistical and Applied Mathematical Sciences Institute (SAMSI). The authors are grateful for the support and inspiring
atmosphere there. The authors would like to thank Prof. P. Degond and Dr. S. Motsch for fruitful and stimulating discussions. The authors
also acknowledge the support by NSF RNMS grant DMS-1107444. The work of A. Chertock was supported in part by the NSF Grant DMS-1115682. The
work of A. Kurganov was supported in part by the NSF Grant DMS-1115718. The work of I. Timofeyev was supported in part by the  NSF Grant
DMS-1109582. I. Timofeyev also acknowledges support from SAMSI as a long-term visitor in Fall 2010 and 2011.
\end{acknowledgment}

\appendix
\section{Central-Upwind Scheme for Systems of Conservation Laws}\label{app}
In this section, we briefly describe a semi-discrete second-order central-upwind scheme used for numerical solution of the
1-D systems \eref{pde} and \eref{pde2}, which can be written as
\begin{equation}
\bm{\rho}_t+\mathbf{F}(\bm{\rho})_x=(Q(\bm{\rho})\bm{\rho}_x)_x,
\label{a1}
\end{equation}
with $\bm{\rho}=(\rho^+,\rho^-)^T$, the flux $\mathbf{F}(\bm{\rho})=\big(f(\rrp)g(\rrm),f(\rrm)g(\rrp)\big)^T$, and the viscosity matrix
$$
Q(\bm{\rho})=\frac{\eps c_0}{2}{\rm diag}\left((1-\rho^-)^2+2\alpha_1\rho^-(1-\rho^-)+\alpha_3(\rho^-)^2,\,
(1-\rho^+)^2+2\alpha_1\rho^+(1-\rho^+)+\alpha_3(\rho^+)^2\right).
$$
For simplicity, we divide the computational domain into finite-volume cells $C_j=[x_{j-\frac{1}{2}},x_{j+\frac{1}{2}}]$ of a uniform size
$\Delta x$ with $x_j=j\Delta x$ and denote by
$$
\overline{\bm{\rho}}_j(t)=\frac{1}{\Delta x}\int\limits_{C_j}\bm{\rho}(x,t)\,dx,
$$
the computed cell averages of the solution $\bm{\rho}$, which we assume to be available at time $t$. The semi-discrete central-upwind 
scheme for \eref{a1} is given by the following system of ODEs for $\overline{\bm{\rho}}_j$:
\begin{equation}
\frac{d\overline{\bm{\rho}}_j(t)}{dt}=-\frac{\mathbf{H_{j-\frac{1}{2}}(t)}-\mathbf{H_{j+\frac{1}{2}}(t)}}{\Delta x}
+\frac{\mathbf{P_{j-\frac{1}{2}}(t)}-\mathbf{P_{j+\frac{1}{2}}(t)}}{\Delta x},
\label{a2}
\end{equation}
where the numerical hyperbolic $\mathbf{H_{j+\frac{1}{2}}(t)}$ and parabolic $\mathbf{P_{j+\frac{1}{2}}(t)}$ fluxes are constructed as
follows (for convenience, we will omit the dependence of the computed quantities on $t$ in the rest of the text).

We first reconstruct a piecewise linear approximation,
\begin{equation}
\widetilde{\bm{\rho}}(x)=\overline{\bm{\rho}}_j+(\bm{\rho}_x)_j(x-x_j),\quad x\in C_j.
\label{r1}
\end{equation}
The numerical derivatives $(\bm{\rho}_x)_j$ should be computed componentwise and are (at least) first-order approximations of 
$\bm{\rho}_x(x_j,t)$, calculated using a nonlinear limiter needed to avoid oscillations in the reconstruction \eref{r1}. In our numerical
experiments, we have used the generalized minmod limiter (see, e.g.,\cite{LN,NT}):
\begin{equation}
(\bm{\rho}_x)_j={\rm minmod}\left(\theta\frac{\overline{\bm{\rho}}_j-\overline{\bm{\rho}}_{j-1}}{\Delta x},\, 
\frac{\overline{\bm{\rho}}_{j+1}-\overline{\bm{\rho}}_{j-1}}{2\Delta x},\,
\theta\frac{\overline{\bm{\rho}}_{j+1}-\overline{\bm{\rho}}_j}{\Delta x}\right),\quad\theta\in[1,2],
\label{minmod}
\end{equation}
where the minmod function is defined as
\begin{equation}
{\rm minmod}(z_1,z_2,...):=\left\{\begin{array}{lc}\!\!\!\min_j\{z_j\}, & ~~\mbox{if} ~~z_j>0 ~~\forall j,\\
\!\!\!\max_j\{z_j\}, & ~~\mbox{if} ~~z_j<0 ~~\forall j,\\ \!\!\!0, & ~\mbox{otherwise},\end{array}\right.
\label{mm}
\end{equation}
and the parameter $\theta$ can be used to control the amount of numerical viscosity present in the resulting scheme (larger
values of $\theta$ correspond to less dissipative but, in general, more oscillatory reconstructions).

The reconstruction \eref{r1} is, in general, discontinuous at the cell interfaces, $x=x_{j+\frac{1}{2}}$, where it has two values, which we
denote by
\begin{equation}
\bm{\rho}_j^{\rm E}:=\overline{\bm{\rho}}_j+\frac{\Delta x}{2}(\bm{\rho}_x)_j,\quad
\bm{\rho}_j^{\rm W}:=\overline{\bm{\rho}}_j-\frac{\Delta x}{2}(\bm{\rho}_x)_j.
\label{a6}
\end{equation}
These discontinuities propagate in time with one-sided local speeds that can be estimated from the eigenvalues $\lambda$ of the
Jacobian matrix \eref{4.8} and are calculated in the following way. We denote by
$$
R=f'(\rrp)g(\rrm)-f'(\rrm)g(\rrp),~~
D=\left[f'(\rrm)g(\rrp)+f'(\rrp)g(\rrm)\right]^2-4f(\rrm)f(\rrp)g'(\rrm)g'(\rrp)
$$
and consider two possible cases:
\begin{itemize}
\item 
If both $D_j^{\rm E}\ge0$ and $D_{j+1}^{\rm W}\ge0$ (hyperbolic regime), then
$$
\begin{aligned}
&a^+_{j+\frac{1}{2}}=\frac{1}{2}\max\limits\left\{R_j^{\rm E}+\sqrt{D_j^{\rm E}},\,R_{j+1}^{\rm W}+\sqrt{D_{j+1}^{\rm W}},\,0\right\},\\
&a^-_{j+\frac{1}{2}}=\frac{1}{2}\min\limits\left\{R_j^{\rm E}-\sqrt{D_j^{\rm E}},\,R_{j+1}^{\rm W}-\sqrt{D_{j+1}^{\rm W}},\,0\right\},
\end{aligned}
$$
\item
If either $D_j^{\rm E}<0$ or $D_{j+1}^{\rm W}<0$ (nonhyperbolic regime), then
$$
\begin{aligned}
&a^+_{j+\frac{1}{2}}=\frac{1}{2}\max\limits\left\{\sqrt{(R_j^{\rm E})^2-D_j^{\rm E}},\,\sqrt{(R_{j+1}^{\rm W})^2-D_{j+1}^{\rm W}}\right\},\\
&a^-_{j+\frac{1}{2}}=-a^+_{j+\frac{1}{2}}.
\end{aligned}
$$
\end{itemize}
The numerical fluxes are then given by 
\begin{equation}
\begin{aligned}
&\mathbf{H_{j+\frac{1}{2}}}=\frac{a^+_{j+\frac{1}{2}}\mathbf{F}(\bm{\rho}_j^{\rm E})
-a^-_{j+\frac{1}{2}}\mathbf{F}(\bm{\rho}_{j+1}^{\rm W})}{a^+_{j+\frac{1}{2}}-a^-_{j+\frac{1}{2}}}
+\frac{a^+_{j+\frac{1}{2}} a^-_{j+\frac{1}{2}}}{a^+_{j+\frac{1}{2}}-a^-_{j+\frac{1}{2}}}
\left[\bm{\rho}_{j+1}^{\rm W}-\bm{\rho}_j^{\rm E}\right],\\[1.1ex]
&\mathbf{P_{j+\frac{1}{2}}}=Q(\bm{\rho}_{j+\frac{1}{2}})\frac{\overline{\bm{\rho}}_{j+1}-\overline{\bm{\rho}}_j}{\Delta x},\quad
\bm{\rho}_{j+\frac{1}{2}}=\frac{\bm{\rho}_j^{\rm E}+\bm{\rho}_{j+1}^{\rm W}}{2}.
\end{aligned}
\label{a8}
\end{equation}
Finally, the resulting semi-discretization \eref{a2}--\eref{a8} is a time-dependent ODE system, which should be numerically integrated using
a stable ODE solver of an appropriate order. In our numerical experiments we have used the third-order strong stability preserving 
Runge-Kutta method (see \cite{GST}).
\begin{remark}
Notice that the choice of one-sided local speeds in the nonhyperbolic regime is ad-hoc. However, it is important to point out that we have
not tried to stabilize the inviscid PDE solution by increasing the amount of numerical viscosity: The solution has been stabilized by adding
nonlinear diffusion terms rigorously derived from the mesoscopic formulation.
\end{remark}

\bibliographystyle{siam}
\bibliography{refs.pt}

\end{document}